\documentclass{amsart}
\usepackage{amsmath,amsthm,amssymb,amsfonts}
\usepackage[colorlinks=true]{hyperref}
\usepackage{color}
\usepackage{graphicx}

\hypersetup{urlcolor=blue, citecolor=blue}%red}

\def\doi#1{   {\href{http://dx.doi.org/#1}
   {{\mdseries\ttfamily DOI}}}}

\newcommand{\R}{\mathbb{R}}

\newcommand{\beeq}{\begin{equation}}\newcommand{\eneq}{\end{equation}}

\newcommand{\supp}{\text{supp}}
\newcommand{\cd}{{\,\cdot\,}}

\newenvironment{prf}{\noindent {\bf Proof.} }{\endprf\par}
\def \endprf{\hfill  {\vrule height6pt width6pt depth0pt}\medskip}

\numberwithin{equation}{section}

\def\O{{\mathcal{O}}}

\def\<{\langle}             \def\>{\rangle}
\def\({\left(}                 \def\){\right)}
\newcommand{\la}{\langle}
\newcommand{\ra}{\rangle}

\newtheorem{thm}{Theorem}[section]

\newtheorem{lem}[thm]{Lemma}

\theoremstyle{remark}
\theoremstyle{definition}
\renewcommand{\S}{{\mathbb{S}}}

\title[Couple wave equations on asymptotically flat backgrounds]
{Global existence for a coupled wave system related to the Strauss conjecture}

\author{Jason Metcalfe}
\address{Department of Mathematics, University of North Carolina,
  Chapel Hill, NC  27599-3250, USA}
\email{metcalfe@email.unc.edu}
\urladdr{http://metcalfe.web.unc.edu}

\author{David Spencer}
\address{University of North Carolina,
  Chapel Hill, NC  27599, USA}
\email{davidspencer6174@gmail.com}

\thanks{
The first author was supported in part by NSF grant DMS-1054289.
The second author was supported in part by a Summer Undergraduate
Research Fellowship (SURF) through the University of North Carolina,
and the results contained herein were developed as a part of his
Undergraduate Honors Thesis.
}
\date{\today}
\dedicatory{} \commby{}

\begin{document}
\bibliographystyle{plain}

\begin{abstract}
A coupled system of semilinear wave equations is considered, and a small data global
existence result related to the Strauss conjecture is
proved.  Previous results have shown that one of the powers may be
reduced below the critical power for the Strauss conjecture provided
the other power sufficiently exceeds such.  The stability of such
results under asymptotically flat perturbations of the space-time
where an integrated local energy decay estimate is available is established.
\end{abstract}

% \keywords{
% Strauss conjecture, Schwarzschild space-time, Kerr space-time, asymptotically flat space-time, weighted Strichartz estimates, localized energy estimates}

% \subjclass[2010]{35L70, 35L15}

\maketitle 
%\tableofcontents

%%%%%%%%%%%%%%%%%%%%%%%%%%%%%%%%%%%%%%%%%%%%%%%%%%
%%%%%%%%%%%%%%%%%%%%%%%%%%%%%%%%%%%%%%%%%%%%%%%%%%
\section{Introduction}

The purpose of this article is to establish global existence for a
coupled system of wave equations, which is related to the Strauss
conjecture, on asymptotically flat space-times that permit a localized
energy estimate.  It is now well-known that nonlinear wave equations
\[\Box u := (\partial_t^2 - \Delta) u = F_p(u),\quad (t,x)\in
\R^+\times \R^n\]
where
\begin{equation}
  \label{Fp}
\sum_{0\le j\le 2} |u|^j |\partial_u^j F_p(u)|\lesssim
|u|^p\quad\text{for $u$ small}
\end{equation}
have global solutions for sufficiently small initial data provided
$p>p_c$ where $p_c>1$ solves
\begin{equation}\label{critical}(n-1)p_c^2 - (n+1)p_c - 2 = 0.\end{equation}
Moreover, blow up is known to occur for $p<p_c$.  These results
originated in \cite{John79} for $n=3$ where $p_c=1+\sqrt{2}$, and
following \cite{Strauss} the problem became known as the Strauss conjecture.  Global existence
in general dimension was eventually established in \cite{GLS97}, \cite{Ta01-2}; see the references
therein for many intermediate results.
Blow up below the critical exponent was proved in \cite{Sideris}.
See also \cite{Schaeffer}, \cite{YZ} for further results at the critical exponent.

In the current work, we shall examine a system of the form
\[\Box u = |v|^p,\quad \Box v = |u|^q.\]
In the flat case, the coupled system was examined in
\cite{dgm}, and it was shown that global existence may be established for
powers in the nonlinearity below the critical exponent provided the
power on the coupled equation exceeds the same.  Indeed, setting 
\begin{equation}
  \label{Cpq}
C(p,q) = \max\Bigl\{\frac{q+2+p^{-1}}{pq-1},
\frac{p+2+q^{-1}}{pq-1}\Bigr\}  - \frac{n-1}{2},
\end{equation}
it was shown that small data global existence holds for $C(p,q)<0$ and
that such fails for $C(p,q)>0$.  Notice that $C(p,p)=0$ corresponds
precisely to \eqref{critical}.  In particular, note that small data
global existence may be established for powers $p<p_c$ provided that
the other power $q$ sufficiently exceeds $p_c$.  In addition to \cite{dgm}, see
\cite{d}, \cite{ds} for treatments of $C(p,q)>0$ and \cite{akt},
\cite{dsm}, \cite{ko_blowup}, and
\cite{GTZ} for analysis of the critical curve $C(p,q)=0$.
Moreover, see the overview \cite{KO} of this and related problems.

Here we
seek to establish the same using techniques that are sufficiently
robust so as to allow background geometries.
Specifically, we shall use a variant of the weighted
Strichartz estimates of \cite{HMSSZ}, \cite{FaWa}, which were further developed in
\cite{LMSTW}, \cite{MWa}, and the localized energy estimate to prove
such global existence.

We shall examine operators of the form
\begin{equation}\label{P}Pu = \partial_\alpha g^{\alpha\beta}\partial_\beta u +
b^\alpha \partial_\alpha u + cu
\end{equation}
on space-times $M$ where $M=\R_+\times \R^3$ or $M=\R_+\times
(\R^3\backslash\mathcal{K})$ where $\mathcal{K}$ has a smooth boundary
and $\mathcal{K}\subset\{x\,:\,|x|<R_0\}$.  Here $g$ is a Lorentzian
metric, and we make the assumption that $g$ can be written as
\begin{equation}
  \label{g}
  g(t,x) = m + g_0(t,r) + g_1(t,x)
\end{equation}
where $m=\text{diag}(-1,1,1,1)$ is the Minkowski metric.  The components
$g_0$ and $g_1$ will represent long-range and short-range
perturbations respectively.  They are asymptotically flat in the sense
that
\begin{equation}
  \label{g_decay}
  \|\partial^\mu_{t,x}
  g_{i,\alpha\beta}\|_{\ell^{i+|\mu|}_1L^\infty_{t,x}} = \O(1),\quad
  i=0,1,\quad |\mu|\le 3.\footnote{Here, as in \cite{MWa}, for a norm
    $A$, we set 
\[\|u\|_{\ell^s_q A} = \Bigl\| 2^{js}
\|\phi_j(x)u(t,x)\|_A\Bigr\|_{\ell^q_{j\ge 0}},\quad
\sum_{j\ge 0} \phi_j^2(x) = 1,\quad \supp\, \phi_j \subset \{\la
x\ra\approx 2^j\}.\]}
\end{equation}
Due to the need to commute with spatial
rotations, the long-range perturbation $g_0$ is assumed to be
spherically symmetric in the sense that the coefficients only
(spatially) depend on $r=|x|$ and 
\begin{equation}
  \label{g_sym}
  g-g_1 = (-1+\tilde{g}_{00}(t,r)) dt^2 + 2
  \tilde{g}_{01}(t,r)\,dt\,dr + (1+\tilde{g}_{11}(t,r)) dr^2 +
  (1+\tilde{g}_{22}(t,r))r^2\, d\omega_{\S^2}^2.
\end{equation}
and, by \eqref{g_decay}, $\|\partial^\mu_{t,x}
\tilde{g}_{\alpha\beta}\|_{\ell^{|\mu|}_1L^\infty_{t,x}}=\O(1)$ for
$|\mu|\le 3$.  
The coefficients of the lower-order perturbations decay are assumed to
decay as follows:
\begin{equation}
  \label{low_order}
  \|\partial^\mu_{t,x} b\|_{\ell_1^{1+|\mu|}L^\infty_{t,x}} +
  \|\partial^\mu_{t,x} c\|_{\ell^{2+|\mu|}_1 L^\infty_{t,x}} =
  \O(1),\quad |\mu|\le 2.
\end{equation}

We shall also assume that the perturbations admit a (weak) localized
energy decay.  More specifically, we assume that there is $R_1$ (with
$R_1>R_0$ in the case that $M=\R_+\times (\R^3\backslash
\mathcal{K})$) so that if $u$ solves $Pu=F$ then
\begin{multline}
  \label{le}
  \|\partial \partial^\mu u\|_{L^\infty_tL^2_x} +
  \|(1-\chi)\partial \partial^\mu u\|_{\ell^{-1/2}_\infty L^2_{t,x}} +
  \|\partial^\mu u\|_{\ell^{-3/2}_\infty L^2_{t,x}}
\\\lesssim \|u(0,\cd)\|_{H^{|\mu|+1}} + \|\partial_t
u(0,\cd)\|_{H^{|\mu|}} + \sum_{|\nu|\le|\mu|} \|\partial^\nu F\|_{L^1_tL^2_x}
\end{multline}
for all $|\mu|\le 2$.  Here $\chi$ is a smooth function that is
identically $1$ on $B_{R_1/2} := \{|x|\le R_1/2\}$ and is supported on
$B_{R_1}$.

On $(1+3)$-dimensional Minkowski space (i.e. when $g_0\equiv g_1\equiv
0$), it is known that
\[\|\partial u\|_{L^\infty_tL^2_x} +
  \|\partial u\|_{\ell^{-1/2}_\infty L^2_{t,x}} +
  \|u\|_{\ell^{-3/2}_\infty L^2_{t,x}}
\lesssim \|\partial u(0,\cd)\|_{L^2} + \|\Box u\|_{L^1_tL^2_x},
\]
which is a stronger version of the $\mu=0$ estimate above.  And as 
the flat d'Alembertian commutes with $\partial_{t,x}$, the higher
order variants readily follow.  Such estimates originated in
\cite{Mo2}.  They follow, e.g., by multiplying $\Box u$ by a multiplier of
the form $\frac{r}{r+2^j}\partial_r u + \frac{n-1}{2}\frac{1}{r+2^j}
u$, integrating over $[0,T]\times \R^3$, and integrating by parts.
See, e.g., \cite{Sterbenz}, \cite{MetSo06}.  And see, e.g., \cite{MTT}
for a more complete history.  These estimates are known to be rather
robust in the asymptotically flat regime.  Even without the cutoff,
they are known to hold for small, possibly time-dependent perturbations
of Minkowski space \cite{MetSo06, MetSo07}, \cite{MetTa07, MetTa09}, \cite{Alinhac}
and for time-independent nontrapping perturbations in the product
manifold setting due to, e.g., \cite{Burq}, \cite{BoHa},
\cite{SoWa10}.  See \cite{MST} for the most general results in the
nontrapping regime.

The presence of trapped rays is a known obstruction to the localized
energy estimate \cite{Ral}, \cite{Sbierski}.  The asymptotic flatness
restricts the possibility of trapped rays to a compact set, and when
the trapping is sufficiently weak, a localized energy estimate where
one, say, cuts off away from the trapping may sometimes be recovered.
Allowing for this is the reason for the cutoff in assumption
\eqref{le}.  Previous results have then verified \eqref{le} in a
number of settings where trapping occurs, including on the
Schwarzschild space-time \cite{BSerr, BS}, \cite{DaRo, DaRo09},
\cite{MMTT}, on Kerr space-times with $a\ll M$ \cite{TT} (see also
\cite{AB09}, \cite{DaRoNew, DaRo08} for some closely related results
and \cite{DaRoSR} for a related result that holds for the full
subextremal range $|a|<M$),
and on certain warped product manifolds that contain degenerate trapping \cite{BCMP}.

We now introduce the specific problem at hand.  With two possibly
different operators $P_1$ and $P_2$ subject to hypotheses
\eqref{g}-\eqref{le}, we examine the coupled system
\begin{equation}
  \label{main_equation}
  \begin{split}
    P_1 u&=F_p(v),\\
u(0,x)&=f_1(x),\\
\partial_tu(0,x)&=g_1(x),
  \end{split}
\qquad\qquad
  \begin{split}
    P_2 v&=F_q(u),\\
v(0,x)&=f_2(x),\\
\partial_tv(0,x)&=g_2(x).
  \end{split}
\end{equation}
Here $F_p$ and $F_q$ are two functions satisfying \eqref{Fp}.

For the system \eqref{main_equation}, we shall establish the following
small data global existence result:

\begin{thm}\label{mainTheorem}
Suppose that $P_1$ and $P_2$ are operators of the form \eqref{P} so
that \eqref{g}, \eqref{g_decay}, \eqref{g_sym}, \eqref{low_order}, and
\eqref{le} hold.  Moreover assume that $2< p, q$ and
$C(p,q)<0$.  Then if $f_1, g_1, f_2, g_2\in C^\infty_c$\footnote{For
  simplicity of exposition, we have taken the data here to be 
compactly supported, but this may be replaced by a condition such as
\cite[(5.3)]{MWa}.} and
\begin{equation}
  \label{data}
  \|(f_1,f_2)\|_{H^3} + \|(g_1,g_2)\|_{H^2} \le \varepsilon
\end{equation}
with $\varepsilon$ sufficiently small, there exists a global solution
$(u,v)$ to \eqref{main_equation}. 
\end{thm}

We note that the techniques to prove Theorem~\ref{mainTheorem} also
work in four spatial dimensions, but as they require $p, q\ge 2$ and
$p_c=2$ in this case, nothing new is gained over \cite{MWa}.

This result follows a number of studies that established various
existence results related to the Strauss conjecture in the presence of
background geometry.  In exterior domains, these included \cite{DMSZ},
\cite{HMSSZ}, \cite{SSW}, and \cite{Yu}.  And on asymptotically flat
backgrounds, see \cite{SoWa10}, \cite{WaYu11}, \cite{LMSTW},
\cite{MWa}, and \cite{Wang}.  See, also, the expository article
\cite{WaYu11p}.  A key component of many of these results is the
weighted Strichartz estimate of \cite{HMSSZ}, \cite{FaWa}.  Here, in
particular, we rely on the variant of that developed in \cite{MWa},
which is based on the local energy estimates of \cite{MetTa07}.

\begin{figure}\caption{Range of allowable indices}\label{fig}
\includegraphics[width=0.7\textwidth]{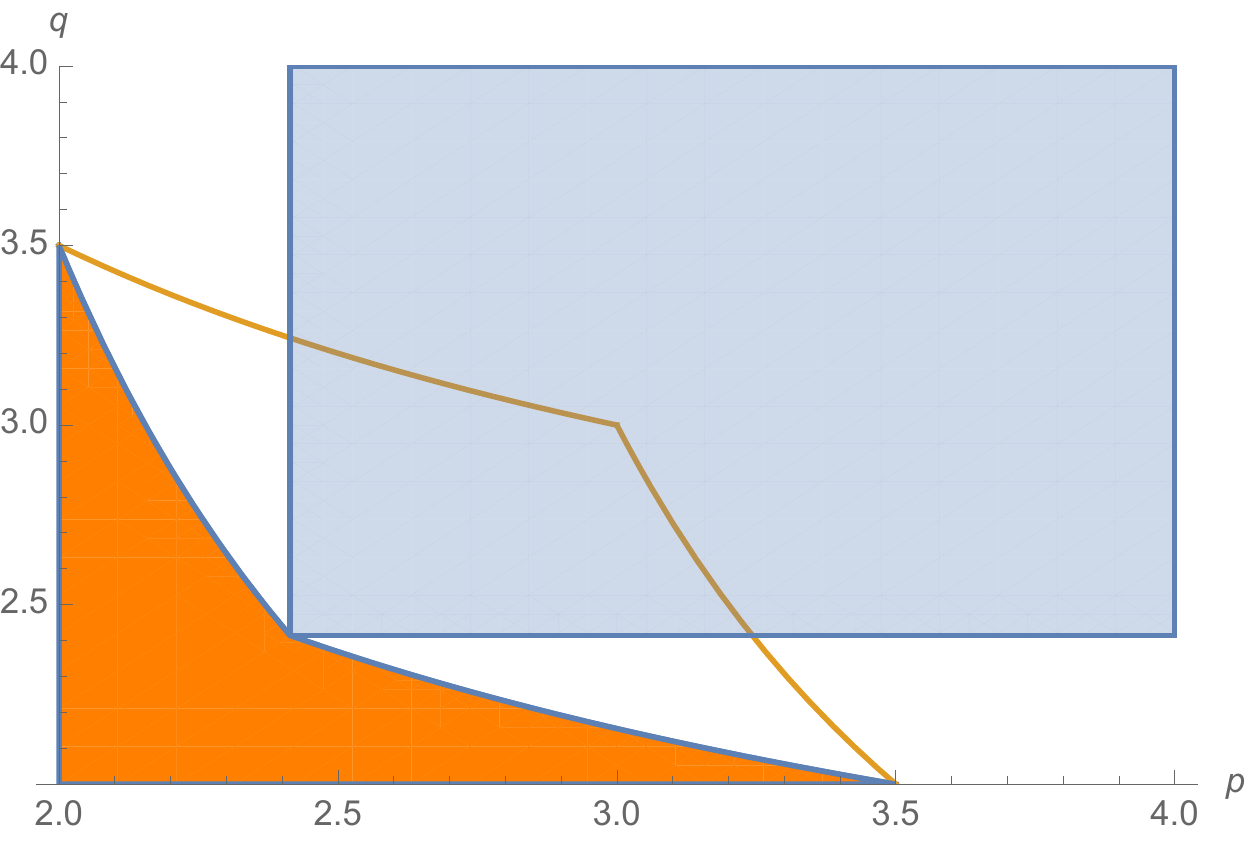}
\end{figure}

Figure~\ref{fig} demonstrates the range of allowable indices.
Theorem~\ref{mainTheorem} allows for any pair of indices $(p,q)$ that
land outside of the shaded region in the lower left corner.  The
techniques of, e.g., \cite{MWa}, however, only trivially apply in the
rectangular region where $p, q > p_c = 1+\sqrt{2}$.  The other curve
will be discussed further in Section~\ref{proof_section}.

The method that we shall employ is similar in spirit to that of
\cite{MWa}, which in turn is based on a number of preceding works.
Near infinity, where the asymptotic flatness allows us to think of the
geometry as a small perturbation of Minkowski space, variants of the
weighted Strichartz estimates of \cite{HMSSZ}, \cite{FaWa} are
employed.  The assumed localized energy estimate handles the remaining
compact region, where the geometry has the most significant role.  It
also allows the analyses done in the two regions to be glued
together.  Such a strategy has become common.  See, e.g.,
\cite{LMSTW}, \cite{MMTT}, \cite{MetTa07, MetTa09}, \cite{MTT},
\cite{Ta13}.

%%%%%%%%%%%%%%%%%%%%%%%%%%%%%%%%%%%%%%%%%%%%%%%%%%%%%%%%%%%%
\section{Main Estimates}

Before we proceed to the main estimates for the linear equation, we
first introduce some notation.  We shall be using a restricted set of
the classical invariant vector fields.  We let $\Omega_{ij}=
x_i\partial_j - x_j\partial_i$ denote the generators of spatial
rotations.  And we let $Y = \{\nabla_x, \Omega\}$ and $Z=\{\partial,
\Omega\}$, where $\partial = (\partial_t,\nabla_x)$ denotes the
space-time gradient.  We also introduce the shorthand $|Z^{\le m} u| =
\sum_{|\mu|\le m} |Z^\mu u|$ for summing over multi-indices of order
$\le m$.  A similar notation where the absolute values are replaced by
a norm shall also be used.  Finally, for the mixed norms that appear
in the weighted Strichartz estimates below, we fix the convention
\[\|f\|_{L^p_t L^q_r L^s_\omega} = \Bigl[\int \Bigl(\int \Bigl[\int
|f(t,r\omega)|^s\,d\omega_{\S^2}\Bigr]^{q/s}
r^2\,dr\Bigr)^{p/q}\,dt\Bigr]^{1/p},\]
with the obvious changes for Lebesgue indices of $\infty$.

The main linear estimate that will be applied near infinity, where
the operators may be viewed as small perturbations of the flat
d'Alembertian, is the following weighted Strichartz estimate:

\begin{thm}[\cite{MWa}]\label{MW_thm}
Suppose that $P$ is an operator of the form \eqref{P} satisfying the
hypotheses \eqref{g}, \eqref{g_decay}, \eqref{g_sym},
\eqref{low_order}, and \eqref{le}.  Suppose that $w(0,\cd)$ and
$\partial_t w(0,\cd)$ are compactly supported and satisfy
\eqref{data}.  Then
there exists $R_2>R_1$ so that for any $R>R_2$ if $\psi_R$ is identically $1$ on
$\{|x|\ge 2R\}$ and vanishes on $\{|x|<R\}$, then we have
\begin{multline}
  \label{wtdStrichartz}
  \|\psi_R Z^{\le 2}
  w\|_{\ell_p^{\frac{3}{2}-\frac{4}{p}-s}L^p_tL^p_rL^2_\omega}
  \le C_1 \varepsilon + C_1\|\psi_R Y^{\le 1} Pw(0,\cd)\|_{\dot{H}^{s-1}}\\+
C_1 \|\psi^{\tilde{p}}_R Z^{\le 2}
  Pw\|_{\ell_1^{-\frac{1}{2}-s}L^1_tL^1_r L^2_\omega} +
 C_1 \|\partial^{\le 2} Pw \|_{L^1_t L^2_r L^2_\omega}
\end{multline}
and
\begin{equation}
  \label{wtdStrichartz0}
  \|\psi_R
  w\|_{\ell_p^{\frac{3}{2}-\frac{4}{p}-s}L^p_tL^p_rL^2_\omega}
  \lesssim \varepsilon + 
 \|\psi^{\tilde{p}}_R 
  Pw\|_{\ell_1^{-\frac{1}{2}-s}L^1_tL^1_r L^2_\omega} +
  \| Pw \|_{L^1_t L^2_r L^2_\omega}
\end{equation}
for any $p\in (2,\infty)$, $s\in (1/2-1/p,1/2)$, and $\tilde{p}>0$. 
\end{thm}

These estimates originated in \cite{HMSSZ} and \cite{FaWa} for the
flat d'Alembertian.  The above version is essentially from \cite{MWa}, which in
particular allows for asymptotically flat operators and does not
necessitate compactly supported initial data though we assume that
here for simplicity.  These estimates follow by interpolating a
variant of the localized energy estimate with a trace theorem on the
sphere.  We note a few minor modifications from the version in
\cite{MWa}: (1) We have allowed a different power of the cutoff
function in the right sides; (2) We have more carefully stated the
requirements on $R$ so that in the sequel we are able to use the same
$R$ for both $u$ and $v$; (3) We have stated separately the case where
no vector fields are applied since this is used when showing that our
iteration is Cauchy.  These all follow from trivial modifications of
the proof of \cite{MWa}.

A key to tying the region near infinity to the remaining compact
region where \eqref{le} is the primary tool is the following weighted
Sobolev inequalities.  These are variants of the original estimates of
\cite{Klainerman} and follow by localizing, applying Sobolev
embeddings on $\R\times \S^2$, and adjusting the volume elements to
match those of $\R^3$ in polar coordinates.  See, e.g., \cite{LMSTW},
\cite{MWa} for proofs. 

\begin{lem}\label{lemmaSob}
  On $\R^3$, for $R\ge 1$, $\beta\in \R$, and $2\le p\le q\le \infty$, we have
  \begin{equation}
    \label{wtdSobolev}
    \|r^\beta u\|_{L^q_r L^\infty_\omega(r\ge R+1)} \lesssim
    \sum_{|\mu|\le 2}
    \|r^{\beta-\frac{2}{p}+\frac{2}{q}} Y^\mu u\|_{L^p_r L^2_\omega
        (r\ge R)}.
  \end{equation}
If $2\le p\le q\le 4$ and $\beta\in \R$, we also have
\begin{equation}
    \label{wtdSobolev2}
    \|r^\beta u\|_{L^q_r L^4_\omega(r\ge R+1)} \lesssim
    \sum_{|\mu|\le 1}
    \|r^{\beta-\frac{2}{p}+\frac{2}{q}} Y^\mu u\|_{L^p_r L^2_\omega
        (r\ge R)}.
  \end{equation}
\end{lem}

%%%%%%%%%%%%%%%%%%%%%%%%%%%%%%%%%%%%%%%%%%%%%%%%%%%%%%%%%%%%
\section{Small data global existence}\label{proof_section}

We now prove Theorem~\ref{mainTheorem}.  We shall apply
\eqref{wtdStrichartz} to $u$ with $(p, s) = \Bigl(q,
\frac{7+4p-3pq}{2-2pq}\Bigr)$ and to $v$ with $(p,s) = \Bigl(p,
\frac{7+4q-3pq}{2-2pq}\Bigr)$.  The requirement from
Theorem~\ref{MW_thm} that $s>\frac{1}{2}-\frac{1}{p}$ then corresponds
to
\[\frac{2+p+q^{-1}}{pq-1}<1,\quad \frac{2+q+p^{-1}}{pq-1}<1\]
respectively.  In the case of $n=3$, this produces the requirement
that $C(p,q)<0$.

We shall assume, without loss of generality, that $(p,q)$ satisfy
\begin{equation}\label{conditions}p(q-2)<3,\quad q(p-2)<3,\end{equation}
which corresponds to the condition that $s<1/2$ in
Theorem~\ref{MW_thm}.  It is these conditions that are represented by
the curve in Figure~\ref{fig}, below which they are satisfied.  In the
unshaded region to the left of the rectangle, if $(p,q)$ satisfy $C(p,q)<0$
but one of the above conditions is violated, one may simply choose any 
$\tilde{q}<q$ so that $C(p,\tilde{q})<0$ and so that both conditions
in \eqref{conditions} hold.  One simply imagines $|u|^q =
|u|^{q-\tilde{q}}|u|^{\tilde{q}}$ and argues as below with the
exponent $\tilde{q}$ replacing $q$.  Simple Sobolev embeddings control
the remaining $q-\tilde{q}$ powers.  In the unshaded region below the
rectangle, one argues similarly by reducing the power of $p$.  In the
shaded rectangle, the methods of \cite{MWa} apply directly and no
further argument is needed.

Let $s_1 = \frac{7+4p-3pq}{2-2pq}$, $s_2 = \frac{7+4q-3pq}{2-2pq}$,
$\alpha_1 = \frac{3}{2}-\frac{4}{q}-s_1$, and $\alpha_2 =
\frac{3}{2}-\frac{4}{p}-s_2$.  We note that the power of the weight in
the right side of \eqref{wtdStrichartz} then satisfies:
\[-\frac{1}{2}-s_1 = p\alpha_2,\quad -\frac{1}{2}-s_2 = q\alpha_1.\]

We solve \eqref{main_equation} via an iteration, and at this point,
the arguments that are used are akin to those of \cite{HMSSZ},
\cite{LMSTW}, and \cite{MWa}.  Setting $u_{-1}
\equiv 0$ and $v_{-1}\equiv 0$, we recursively define $u_j, v_j$,
$j\ge 0$ to
solve
\begin{equation}
  \label{iteration}
  \begin{split}
    P_1 u_j&=F_p(v_{j-1}),\\
u_j(0,x)&=f_1(x),\\
\partial_tu_j(0,x)&=g_1(x),
  \end{split}
\qquad\qquad
  \begin{split}
    P_2 v_j&=F_q(u_{j-1}),\\
v_j(0,x)&=f_2(x),\\
\partial_tv_j(0,x)&=g_2(x).
  \end{split}
\end{equation}
We introduce the quantity
\begin{multline}
  \label{norm}
  M_k(u,v) = \|\psi_R Z^{\le k} u\|_{\ell^{\alpha_1}_q L^q_tL^q_r L^2_\omega} +
  \|\psi_R Z^{\le k} v\|_{\ell^{\alpha_2}_pL^p_tL^p_r L^2_\omega} \\+
  \|\partial^{\le k} (u,v)\|_{\ell_\infty^{-3/2}L^2_tL^2_rL^2_\omega}
  + \|\partial^{\le k} \partial (u,v)\|_{L^\infty_tL^2_rL^2_\omega}.
\end{multline}

Our first goal is to inductively show that $M_2(u_j,v_j) \le
4C_2\varepsilon$ for some uniform constant $C_2$.  We first note that
from \eqref{data} (and the assumption that the data are compactly supported) we easily obtain
\[\|\psi_R Y^{\le 1} P_1u_j(0,\cd)\|_{\dot{H}^{s_1-1}} + \|\psi_R Y^{\le
  1} P_2 v_j(0,\cd)\|_{\dot{H}^{s_2-1}} \le C\varepsilon.\]
And thus, for $C_2$ chosen large enough, by \eqref{wtdStrichartz},
\begin{multline}\label{StrichartzApp}
  \|\psi_R Z^{\le 2} u_j\|_{\ell^{\alpha_1}_q L^q_tL^q_r L^2_\omega} +
  \|\psi_R Z^{\le 2} v_j\|_{\ell^{\alpha_2}_pL^p_tL^p_r L^2_\omega} \le
  C_2\varepsilon \\+ C \|\psi^p_R Z^{\le 2}
  F_p(v_{j-1})\|_{\ell^{p\alpha_2}_1 L^1_tL^1_rL^2_\omega} + C
  \|\psi^q_R Z^{\le 2} F_q(u_{j-1})\|_{\ell^{q\alpha_1}_1
    L^1_tL^1_rL^2_\omega} \\+ C \|\partial^{\le 2}
  F_p(v_{j-1})\|_{L^1_tL^2_rL^2_\omega} + C \|\partial^{\le 2} F_q(u_{j-1})\|_{L^1_tL^2_rL^2_\omega}.
\end{multline}
And from \eqref{data} and \eqref{le}, we have
\begin{multline}
  \label{leApp}
  \|\partial^{\le 2}(u_j,v_j)\|_{\ell^{-3/2} L^2_tL^2_rL^2_\omega} +
  \|\partial^{\le 2} \partial(u_j,v_j)\|_{L^\infty_tL^2_rL^2_\omega} \le
  C_2\varepsilon \\+ C\|\partial^{\le 2}
  F_p(v_{j-1})\|_{L^1_tL^2_rL^2_\omega} + C\|\partial^{\le 2} F_q(u_{j-1})\|_{L^1_tL^2_rL^2_\omega}.
\end{multline}
From these, we first notice that $M_2(u_0,v_0)\le 2C_2\varepsilon$,
which provides the base case for the induction.

We then assume that $M_2(u_{j-1}, v_{j-1}) \le 4C_2\varepsilon$ and
show that $M_2(u_j, v_j)\le 4C_2\varepsilon$.  We first notice that
\eqref{Fp} gives
\begin{equation}\label{prodRule}|Z^{\le 2} F_p(u)|\lesssim |u|^{p-1} |Z^{\le 2} u| + |u|^{p-2}
|Z^{\le 1} u|^2.\end{equation}
By the Sobolev embeddings $H^2_\omega\subset L^\infty_\omega$ and
$H^1_\omega\subset L^4_\omega$ on $\S^2$ and H\"older's inequality,
this gives
\begin{align*}
  \|Z^{\le 2} F_p(v_{j-1})\|_{L^2_\omega}&\lesssim
  \|v_{j-1}\|_{L^\infty_\omega}^{p-1} \|Z^{\le 2}
  v_{j-1}\|_{L^2_\omega} + \|v_{j-1}\|_{L^\infty_\omega}^{p-2}
\|Z^{\le 1} v_{j-1}\|_{L^4_\omega}^2\\
&\lesssim \|Z^{\le 2} v_{j-1}\|^p_{L^2_\omega}.
\end{align*}
It then immediately follows that
\begin{equation}
  \label{I}
  \|\psi^p_R Z^{\le 2} F_p(v_{j-1})\|_{\ell^{p\alpha_2}_1
    L^1_tL^1_rL^2_\omega} \lesssim \|\psi_R Z^{\le
    2}v_{j-1}\|^p_{\ell^{\alpha_2}_p L^p_tL^p_rL^2_\omega} \lesssim
  (M_2(u_{j-1},v_{j-1}))^p\lesssim \varepsilon^p.
\end{equation}
The last inequality results from the inductive hypothesis.  A similar
argument shows that
\begin{equation}
  \label{II}
  \|\psi^q_R Z^{\le 2} F_q(u_{j-1})\|_{\ell^{q\alpha_1}_1
    L^1_tL^1_rL^2_\omega} \lesssim \varepsilon^q.
\end{equation}

To finish the proof of the boundedness of $M_2(u_j,v_j)$, it remains to
examine the $L^1_tL^2_rL^2_\omega$ pieces in the right sides of
\eqref{StrichartzApp} and \eqref{leApp}.  The analyses will be done
separately outside of a ball of radius $2R+1$, where the resulting
terms will be compared to the weighted Strichartz portions of
$M_2(u_{j-1},v_{j-1})$, and inside the remaining compact set where the
localized energy portions will be used.

We begin with the former using the weighted Sobolev inequalities of
Lemma~\ref{lemmaSob}.  Starting again at \eqref{prodRule}, we have
\begin{multline}\label{energyExterior}
  \|Z^{\le 2}F_p(v_{j-1})\|_{L^1_tL^2_{r\ge 2R+1} L^2_\omega} \lesssim  
  \|r^{-\frac{\alpha_2}{p-1}} v_{j-1}\|^{p-1}_{L^p_t
    L^{\frac{2p(p-1)}{p-2}}_{r\ge 2R+1}L^\infty_\omega} \|r^{\alpha_2} \psi_R
  Z^{\le 2} v_{j-1}\|_{L^p_tL^p_r L^2_\omega}
\\+\Bigl\|r^{\frac{2}{p-2}\bigl(-\alpha_2 - \frac{2}{p}+\frac{1}{2}\bigr)}
v_{j-1}\Bigr\|^{p-2}_{L^p_t L^\infty_{r\ge 2R+1} L^\infty_\omega}
\|r^{\alpha_2 + \frac{2}{p}-\frac{1}{2}} Z^{\le 1} v_{j-1}\|^2_{L^p_t
  L^4_{r\ge 2R+1} L^4_\omega}.
\end{multline}
Since \eqref{conditions} guarantees that $s_1<1/2$, it follows that
$p\alpha_2 = -\frac{1}{2}-s_1 \ge -1$, which is equivalent to
\[-\frac{\alpha_2}{p-1} -\frac{2}{p} + \frac{p-2}{p(p-1)}\le
\alpha_2.\]
Thus, \eqref{wtdSobolev} yields
\[\|r^{-\frac{\alpha_2}{p-1}} v_{j-1}\|_{L^p_t
    L^{\frac{2p(p-1)}{p-2}}_{r\ge 2R+1}L^\infty_\omega} 
\lesssim \|r^{\alpha_2} Z^{\le 2} v_{j-1}\|_{L^p_t L^p_{r\ge
    2R}L^2_\omega}.\]
The estimate \eqref{wtdSobolev2} directly applies to yield
\[\|r^{\alpha_2 + \frac{2}{p}-\frac{1}{2}} Z^{\le 1} v_{j-1}\|_{L^p_t
  L^4_{r\ge 2R+1} L^4_\omega}\lesssim \|r^{\alpha_2} Z^{\le 2}
v_{j-1}\|_{L^p_tL^p_{r\ge 2R}L^2_\omega}.\]
As above, $p\alpha_2 \ge -1$, which implies that
\[\frac{2}{p-2}\Bigl(-\alpha_2
-\frac{2}{p}+\frac{1}{2}\Bigr)-\frac{2}{p}\le \alpha_2.\]  Hence
\eqref{wtdSobolev} gives
\[\Bigl\|r^{\frac{2}{p-2}\bigl(-\alpha_2 - \frac{2}{p}+\frac{1}{2}\bigr)}
v_{j-1}\Bigr\|_{L^p_t L^\infty_{r\ge 2R+1} L^\infty_\omega} \lesssim
\|r^{\alpha_2} Z^{\le 2} v_{j-1}\|_{L^p_tL^p_{r\ge 2R}L^2_\omega}.\]
Plugging each of these bounds into \eqref{energyExterior} and using
that $\psi_R$ is identically 1 on $r\ge 2R$, it follows that
\begin{equation}\label{III} \|Z^{\le 2}F_p(v_{j-1})\|_{L^1_tL^2_{r\ge 2R+1} L^2_\omega}
\lesssim \|\psi_R
  Z^{\le 2} v_{j-1}\|^p_{\ell^{\alpha_2}_p L^p_tL^p_r
    L^2_\omega} \lesssim (M_2(u_{j-1},v_{j-1}))^p \lesssim
  \varepsilon^p.\end{equation}
And an analogous argument in the symmetric variable $q$ shows that
\begin{equation}
  \label{IV}
  \|Z^{\le 2}F_q(u_{j-1})\|_{L^1_tL^2_{r\ge 2R+1} L^2_\omega}
\lesssim 
  \varepsilon^q,
\end{equation}
which leaves the analysis of the $L^1_tL^2_rL^2_\omega$ pieces over the
region $r\le 2R+1$ where the coefficients of $Z$ are bounded.

We again start with \eqref{prodRule} and apply H\"older's inequality
to obtain
\begin{multline}\label{energyInterior}
  \|Z^{\le 2} F_p(v_{j-1})\|_{L^1_tL^2_{r\le 2R+1} L^2_\omega}
\\  \lesssim \|v_{j-1}\|^{p-2}_{L^\infty_t L^\infty_r L^\infty_\omega}
  \|v_{j-1}\|_{L^2_t L^\infty_{r\le 2R+1} L^\infty_\omega} \|\partial^{\le 2}
  v_{j-1}\|_{L^2_t L^2_{r\le 2R+1}L^2_\omega}
\\+\|v_{j-1}\|^{p-2}_{L^\infty_t L^\infty_r L^\infty_\omega} \|\partial^{\le
  1} v_{j-1}\|^2_{L^2_t L^4_{r\le 2R+1} L^4_\omega}.
\end{multline}
By Sobolev embeddings, we have
\[\|v_{j-1}\|_{L^\infty_t L^\infty_r L^\infty_\omega} \lesssim
\|\partial^{\le 1} v_{j-1}\|_{L^\infty_t L^6_r L^6_\omega}\lesssim
\|\partial^{\le 1} \partial v_{j-1}\|_{L^\infty_t L^2_rL^2_\omega}.\]
Similarly, Sobolev embeddings (with a localizing factor) give
\[\|v_{j-1}\|_{L^2_tL^\infty_{r\le 2R+1} L^\infty_\omega} \lesssim
\|\partial^{\le 2} v_{j-1}\|_{L^2_t L^2_{r\le 2R+2} L^2_\omega}
\lesssim \|\partial^{\le 2} v_{j-1}\|_{\ell^{-3/2}_\infty L^2_t L^2_r
  L^2_\omega}\]
and
\[\|\partial^{\le 1} v_{j-1}\|_{L^2_t L^4_{r\le 2R+1} L^4_\omega} \lesssim
\|\partial^{\le 2} v_{j-1}\|_{L^2_t L^2_{r\le 2R+2} L^2_\omega}
\lesssim \|\partial^{\le 2} v_{j-1}\|_{\ell^{-3/2}_\infty L^2_t L^2_r
  L^2_\omega}.\]
These bounds in \eqref{energyInterior} show that
\begin{equation}
  \label{V}\begin{split}
  \|Z^{\le 2}F_p(v_{j-1})\|_{L^1_t L^2_{r\le 2R+1}L^2_\omega} &\lesssim
  \|\partial^{\le 1}\partial v_{j-1}\|^{p-2}_{L^\infty_t L^2_r
    L^2_\omega} \|\partial^{\le 2} v_{j-1}\|^2_{\ell^{-3/2}_\infty
    L^2_tL^2_rL^2_\omega} \\&\lesssim (M_2(u_{j-1},v_{j-1}))^p \lesssim \varepsilon^p.
\end{split}
\end{equation}
And an analogous argument gives
\begin{equation}
  \label{VI}
   \|Z^{\le 2}F_q(u_{j-1})\|_{L^1_t L^2_{r\le 2R+1}L^2_\omega}
   \lesssim \varepsilon^q.
\end{equation}

Using \eqref{I}, \eqref{II}, \eqref{III}, \eqref{IV}, \eqref{V}, and
\eqref{VI} in \eqref{StrichartzApp} and \eqref{leApp}, we obtain
\[M_2(u_j,v_j) \le 2C_2\varepsilon + C_3 \varepsilon^{\min(p,q)}\]
for some constant $C_3$ that is independent of $j$.  Since $p, q\ge
2$, if $\varepsilon$ is sufficiently small, the desired bound
$M_2(u_j,v_j)\le 4C_2\varepsilon$ follows.

It remains to show that the sequence converges.  To do so, we shall
show
\begin{equation}\label{CauchyGoal}
M_0(u_j-u_{j-1}, v_j-v_{j-1}) \le \frac{1}{2} M_0(u_{j-1}-u_{j-2}, v_{j-1}-v_{j-2}),\end{equation}
and this will complete the proof.  Applying \eqref{wtdStrichartz0} and
\eqref{le}, we have
\begin{multline}\label{StrichartzApp2}
  \|\psi_R (u_j-u_{j-1})\|_{\ell^{\alpha_1}_q L^q_tL^q_r L^2_\omega} +
  \|\psi_R (v_j-v_{j-1})\|_{\ell^{\alpha_2}_pL^p_tL^p_r
    L^2_\omega} \\\lesssim 
\|\psi^p_R 
  (F_p(v_{j-1})-F_p(v_{j-2})\|_{\ell^{p\alpha_2}_1 L^1_tL^1_rL^2_\omega} + 
  \|\psi^q_R (F_q(u_{j-1})-F_q(u_{j-2}))\|_{\ell^{q\alpha_1}_1
    L^1_tL^1_rL^2_\omega} \\+  \|
  F_p(v_{j-1}) - F_p(v_{j-2})\|_{L^1_tL^2_rL^2_\omega} +  \|F_q(u_{j-1})-F_q(u_{j-2})\|_{L^1_tL^2_rL^2_\omega},
\end{multline}
and
\begin{multline}
  \label{leApp2}
  \|(u_j-u_{j-1},v_j-v_{j-1})\|_{\ell^{-3/2} L^2_tL^2_rL^2_\omega} +
  \|\partial(u_j-u_{j-1},v_j-v_{j-1})\|_{L^\infty_tL^2_rL^2_\omega} \\\lesssim
  \|
  F_p(v_{j-1})-F_p(v_{j-2})\|_{L^1_tL^2_rL^2_\omega} + \|F_q(u_{j-1})-F_q(u_{j-2})\|_{L^1_tL^2_rL^2_\omega}.
\end{multline}

Noting that
\begin{equation}
  \label{taylor}
  |F_p(v_{j-1})-F_p(v_{j-2})|\lesssim (|v_{j-1}|^{p-1} +
  |v_{j-2}|^{p-1}) |v_{j-1}-v_{j-2}|,
\end{equation}
we can quickly observe that
\begin{equation}
  \label{Ia}
\begin{split}  \|\psi^p_R&(F_p(v_{j-1}) -F_p(v_{j-2})\|_{\ell^{p\alpha_2}_1
    L^1_tL^1_rL^2_\omega} \\&\lesssim
  \Bigl(\|\psi_R v_{j-1}\|^{p-1}_{\ell^{\alpha_2}_p L^p_tL^p_rL^\infty_\omega} +
  \|\psi_R v_{j-2}\|^{p-1}_{\ell^{\alpha_2}_p
    L^p_tL^p_rL^\infty_\omega}\Bigr)
  \|\psi_R(v_{j-1}-v_{j-2})\|_{\ell^{\alpha_2}_p
    L^p_tL^p_rL^2_\omega}\\
&\lesssim \Bigl(\|\psi_R Z^{\le 2} v_{j-1}\|^{p-1}_{\ell^{\alpha_2}_p L^p_tL^p_rL^2_\omega} +
  \|\psi_R Z^{\le 2}v_{j-2}\|^{p-1}_{\ell^{\alpha_2}_p
    L^p_tL^p_rL^2_\omega}\Bigr)\\&\qquad\qquad\qquad\qquad\qquad\qquad\qquad\qquad\qquad\times
  \|\psi_R(v_{j-1}-v_{j-2})\|_{\ell^{\alpha_2}_p
    L^p_tL^p_rL^2_\omega}\\
&\lesssim (4C_2\varepsilon)^{p-1} M_0(u_{j-1}-u_{j-2},v_{j-1}-v_{j-2}).
\end{split}
\end{equation}
A similar argument yields
\begin{equation}
  \label{IIa}
  \|\psi^q_R(F_q(u_{j-1}) -F_q(u_{j-2})\|_{\ell^{q\alpha_1}_1
    L^1_tL^1_rL^2_\omega} 
\lesssim (4C_2\varepsilon)^{q-1} M_0(u_{j-1}-u_{j-2},v_{j-1}-v_{j-2}).
\end{equation}

We also start at \eqref{taylor} to control the $L^1_tL^2_rL^2_\omega$
terms.  There we get, using \eqref{wtdSobolev} as above,
\begin{equation}
  \label{IIIa}
  \begin{split}
    \|&F_p(v_{j-1})-F_p(v_{j-2})\|_{L^1_tL^2_{r\ge 2R+1} L^2_\omega} \\&\lesssim
    (\|r^{-\frac{\alpha_2}{p-1}} v_{j-1}\|^{p-1}_{L^p_t
      L^{\frac{2p(p-1)}{p-2}}_{r\ge 2R+1} L^\infty_\omega}
  + \|r^{-\frac{\alpha_2}{p-1}} v_{j-2}\|^{p-1}_{L^p_t
      L^{\frac{2p(p-1)}{p-2}}_{r\ge 2R+1} L^\infty_\omega})
    \|r^{\alpha_2} \psi_R (v_{j-1}-v_{j-2})\|_{L^p_tL^p_rL^2_\omega}\\
&\lesssim [ (M_2(u_{j-1},v_{j-1}))^{p-1} +
(M_2(u_{j-2},v_{j-2}))^{p-1}] M_0(u_{j-1}-u_{j-2}, v_{j-1}-v_{j-2})\\
&\lesssim \varepsilon^{p-1} M_0(u_{j-1}-u_{j-2}, v_{j-1}-v_{j-2}).
  \end{split}
\end{equation}
And similarly
\begin{equation}
  \label{IVa}
    \|F_q(u_{j-1})-F_q(u_{j-2})\|_{L^1_tL^2_{r\ge 2R+1} L^2_\omega} 
\lesssim \varepsilon^{q-1} M_0(u_{j-1}-u_{j-2}, v_{j-1}-v_{j-2}).
\end{equation}
Finally, using Sobolev embeddings as in \eqref{V},
\begin{equation}
  \label{Va}
  \begin{split}
    \|&F_p(v_{j-1})-F_p(v_{j-2})\|_{L^1_tL^2_{r\le 2R+1} L^2_\omega} \\&\lesssim
 [ \|v_{j-1}\|^{p-2}_{L^\infty_tL^\infty_rL^\infty_\omega}
   \|v_{j-1}\|_{L^2_tL^\infty_{r\le 2R+1}L^\infty_\omega} +\|v_{j-2}\|^{p-2}_{L^\infty_tL^\infty_rL^\infty_\omega}
   \|v_{j-2}\|_{L^2_tL^\infty_{r\le 2R+1}L^\infty_\omega}
   ]\\&\qquad\qquad\qquad\qquad\qquad\qquad\qquad\times 
   \|v_{j-1}-v_{j-2}\|_{L^2_t L^2_{r\le 2R+1} L^2_\omega}
\\
&\lesssim [ (M_2(u_{j-1},v_{j-1}))^{p-1} +
(M_2(u_{j-2},v_{j-2}))^{p-1}] M_0(u_{j-1}-u_{j-2}, v_{j-1}-v_{j-2})\\
&\lesssim \varepsilon^{p-1} M_0(u_{j-1}-u_{j-2}, v_{j-1}-v_{j-2}),
  \end{split}
\end{equation}
and by the same procedures
\begin{equation}
  \label{VIa}
   \|F_q(u_{j-1})-F_q(u_{j-2})\|_{L^1_tL^2_{r\le 2R+1} L^2_\omega} 
\lesssim \varepsilon^{q-1} M_0(u_{j-1}-u_{j-2}, v_{j-1}-v_{j-2}).
\end{equation}

Plugging \eqref{Ia}-\eqref{VIa} into \eqref{StrichartzApp2} and
\eqref{leApp2} immediately yields \eqref{CauchyGoal} provided that
$\varepsilon$ is sufficiently small, and this completes the proof.

%%%%%%%%%%%%%%%%%%%%%%%%%%%%%%%%%%%%%%%%%%%%%%%%%%%%%%%%%%%%% 
\bibliography{strauss_bib}

\end{document}